\documentclass[12pt,oneside]{article}
\usepackage{amsmath,amssymb,amsfonts,amsthm}
\usepackage{hyperref}
\hypersetup{
	colorlinks=true,
	linkcolor=blue,
	filecolor=magenta,
	urlcolor=cyan,
	citecolor=red}
\newcommand{\T}{{\cal T}}

\newcommand{\Real}{\mathbb R}

\newcommand{\set}[1]{\left\{#1\right\}}

\newcommand{\To}{\longrightarrow}

\newcommand {\cp}{\mathfrak{X}(\pi )}

\def\pa{\partial}

\newtheorem{thm}{Theorem}[section]
\newtheorem{cor}[thm]{Corollary}
\newtheorem{lem}[thm]{Lemma}
\newtheorem{prop}[thm]{Proposition}
\newtheorem{defn}[thm]{Definition}
\newtheorem{example}{Example}
\newtheorem{rem}[thm]{Remark}

\numberwithin{equation}{section}

\begin{document}
\title{{{\bf On Generalized Matsumoto Metrics with a Special $\pi$-form }}}
\author{\bf{ A. Soleiman$^1$ and Ebtsam H. Taha$^2$ }}
\date{}
\maketitle                     
\vspace{-1.0cm}
\begin{center}
{$^{1}$ Department of Mathematics, Faculty of Science,  Jouf University,  Skaka,  Kingdom of Saudia Arabia
\vspace{0.2cm}
\\ $^{2}$ Department of Mathematics, Faculty of Science, Cairo University, Giza,  Egypt}
\end{center}
\vspace{-0.5cm}
\begin{center}
E-mails:$^{1}$ asoliman@ju.edu.sa,  amrsoleiman@yahoo.com\\
{\hspace{1.8cm}}$^{2}$ ebtsam@sci.cu.edu.eg, ebtsam.h.taha@hotmail.com
\end{center}

\begin{center}
\textbf{\textit{Dedicated to the memory of Professor Nabil L. Youssef}}
\end{center}

\vspace{0.7cm} \maketitle
\smallskip

\noindent{\bf Abstract.} We explore a generalization of  Matsumoto metric intrinsically.  Given a Finsler manifold $(M,F)$ which admits a concurrent $\pi$-vector field $\overline{\varphi}$, we consider the change $\widehat{F}(x,y)=\frac {F^2 (x,y)} {F(x,y)-\Phi(x,y)}$, where $\Phi$ is  the associated concurrent $\pi$-form with $F(x,y) > \Phi(x,y)$ for all $(x,y) \in \T M$.  We find the condition under which the generalized $\phi$-Matsumoto metric $\widehat{F}$ is a Finsler metric. Moreover, the relations between the associated Finslerian  geometric objects of $\widehat{F}$ and $F$ are obtained, namely,  the relations between angular metric tensors, metric tensors, Cartan torsions, geodesic sprays,   Barthel connections (along with its curvature) and Berwald connections. Further, we prove that the Finsler metrics $F$ and $\widehat{F}$ can never be projectively related. Also, a  condition for the $\pi$-vector field $\overline{\varphi}$ to be concurrent with respect to $\widehat{F}$ is acquired.  Moreover,  an example  of a rational Finsler metric admitting a concurrent $\pi$-vector field together with  the associated change $\widehat{F}$ is provided. Finally,  we find the conditions that preserve  the almost rationality property of a Finsler metric $F$ under the $\phi$-Matsumoto change.



\vspace{30truept}\centerline{\Large\bf{Introduction}}\vspace{12pt}
\par

The Randers metric is defined by $F=\alpha+B$, with $\alpha$ is defined by a Riemannian metric  and $B$ is a $1$-form on the manifold $M$. It named after G.  Randers, 1941,  who introduced a  simple Finslerian metric by this  change.   It is used to construct a generalized field theory that
would comprise gravity and electromagnetism.  Randers metrics are important as they represent solutions to  Zermelo’s  navigation problem in the case of a weak wind (see, e.g. \cite{Erasmo} for a wide perspective of Randers or the more general notion of a wind Finsler metric).  
In  1974,  M.  Matsumoto studied Randers spaces in a more general setting,  by assuming that $\alpha$ is a Finsler metric \cite{r2.8}.  Later,   generalized Randers metrics were investigated in \cite{Tamim-Youssef}.  Another interesting metric is the Matsumoto (Matsumoto's slope of a mountain) metric defined by $F= \frac{\alpha^2}{\alpha-B}$, where  $B$  was originally taken to be induced by earth's gravity.  It was introduced by Matsumoto as a realization of the idea of a slope measure of a mountain with respect to a time measure \cite{ [6],[3]}.  A further local study of Matsumoto metric and its generalizations can be found in \cite{h-Matsumoto change, Matsumoto Change of Finsler Metric,  Tayebi1, Tayebi2, Tayebi3}. On the other hand,  the existence of a concurrent vector field on Finsler spaces has been studied firstly by Tachibana \cite{Tachibana}.  The existence of a concurrent vector field is a very rigid property as, for example, a $3$-dimensional Finsler manifold, a Finsler surface,  a Landsberg space, a $C$-reducible Finsler space admitting a concurrent vector field is Riemannian  \cite{MatEguchi}. 
Also, a generalization of  a concurrent vector field,  which is called  a semi-concurrent vector field, has been investigated in~\cite{Semi-concurrent}.

\medskip

In this paper, we provide an intrinsic investigation of what we called a \textit{generalized $\phi$-Matsumoto metric}. Our intrinsic formulation and index-free proofs in the first three sections  give rise to  simple compact results that  hold globally on the manifold.  More precisely,  in \S 2, we provide a coordinate-free study of a generalized Matsumoto metric with a special  $\pi$-form.  By a generalized Matsumoto metric we mean  the change of a Finsler metric $F$ (not necessarily Riemannian) by a $1$-form $B$,  defined by,  $\widehat{F}=\frac {F^2} {F-B}$, with $F(x,y) > \Phi(x,y)$ for all $(x,y) \in \T M$.  We consider a Finsler manifold $(M,F)$ that admits a concurrent $\pi$-vector field $\overline{\varphi}$ and find the corresponding $\pi$-form $\Phi$. Thus the associated  $1$-form ${\Phi}(x,y)$ is used to define the $\phi$-Matsumoto change \eqref{change}.
We  analyse intrinsically  some of the geometric objects associated with $\widehat{F}$,  namely, the supporting form $\widehat{\ell}$, the angular metric tensor $\widehat{\hbar}$, the Finsler metric tensor $\widehat{g}$
 and the Cartan torsion $\widehat{\mathbf{T}}$. Also, we characterize the non-degeneracy property of the metric tensor $\widehat{g}$, that is,  $\widehat{g}$ is non-degenerate if and only if $F(1+2g(\overline{\varphi},\overline{\varphi}))-3\Phi\neq0$.  Then, in \S 3,    the related geodesic sprays corresponding to this change as well as the relation between the two Barthel connections $\Gamma$ and $\widehat{\Gamma}$ are obtained.  As a consequence, we prove that two sprays are projectively related if and only if  $\overline{\varphi}=0$, which contradicts the assumption of $\pi$-concurrent vector field.   Then, the generalized Matsumoto change of  Berwald connection  implies that the $\pi$-vector field $\overline{\varphi}$ is never concurrent with respect to the Finsler metric $\widehat{F}$.  Finally, in \S 4,  an example of a rational Finsler metric $F$ that admits a $\pi$-concurrent vector field is given and its  $\phi$-Matsumoto change is obtained.  Finally,  we tried to answer the question: if $F$ is  an almost rational Finsler metric,  under what conditions $\widehat{F}$ is an almost rational Finsler metric?
\section{Preliminaries of Finsler geometry}
\par Let $M$ be an $n$-dimensional para-compact smooth manifold and  $\pi: T M\longrightarrow M$ its  tangent bundle. The vertical subbundle $V(TM)$ is defined to be $ \ker (d\pi)$.  Let $\T M:=TM/\set{0}$ be the slit tangent bundle, $\dot{\pi}: \T M\longrightarrow M$ and  the pullback bundle of the tangent bundle be denoted by $\dot{\pi}^{-1}(T M)$ over $\T M$. It is called also the Finsler bundle \cite{2014}.  Further,  $\mathfrak{F}(\T M)$ denotes the algebra of smooth functions on $\T M$ and $\cp$ the $\mathfrak{F}(\T M)$-module of differentiable sections of $\dot{\pi}^{-1}(\T M)$. The elements of $\mathfrak{X}(\pi )$ will be called $\pi$-vector fields and denoted by barred letters $\overline{X}$. \\
 \vspace*{-0.4cm}
\par
We have the short exact sequence \cite{r21, r94a}
$$0\longrightarrow
 \dot{\pi}^{-1}(\T M)\stackrel{\gamma}\longrightarrow T\T M \stackrel{\rho}\longrightarrow
\dot{\pi}^{-1}(\T M)\longrightarrow 0 ,\vspace{-0.1cm}$$ where $\gamma$ is the natural injection ($\gamma: \dot{\pi}^{-1}(\T M)\rightarrow V(\T M)$ is an isomorphism) and $\rho := (\dot{\pi}, d \pi)$. The vertical endomorphism  $J$  is the map $J: T\T M \longrightarrow T\T M$ defined by $J=\gamma \circ \rho$. For all  $ f \in \mathfrak{F}(\T M), \, W \in \mathfrak{X}(\T M)$, $J$ satisfies:
 \begin{equation}\label{eq. 1}
[f W, J]=f[W,J]+df \wedge i_{W}J-d_{J}f\otimes W,
\end{equation}
where $d f$ is the exterior derivative of $f$,  the derivative $d_{J}:=[i_{J} ,d]= i_{J} \circ d - d \circ i_{J}$ and
 $i_{W}$ is the interior product with respect to $ W$ defined by $i_{W} f = 0, \,\, i_{W}J= J(W)$. Moreover, we have:
\begin{equation}\label{J properties}
 i_{ \mathcal{C}} \,J=0 \quad \text{    and    } \,\, \, [\mathcal{C}, J] = -J,
 \end{equation}
where the vector field defined by $\mathcal{C}:=\gamma\, \overline{\eta}$ is called Liouville vector field and  $\overline{\eta}(u)=(u,u)$ for all $u$ in $\T M .$ Moreover,  a spray on $M$ is a smooth vector field $G$ on $\T M$ such that  $J G = \mathcal{C}$ and $[\mathcal{C},G]= G$.
\bigskip
\par
For a linear connection  $D$ on $\dot{\pi}^{-1}(\T M)$, we have    $K:T\T M \longrightarrow
\dot{\pi}^{-1}(\T M)$ which is defined by $K(W) =D_W \overline{\eta}$.  Thereby,  the horizontal space at $u \in TM$ is $H_u (\T M):= \{ W \in T_u
(\T M) : K(W)=0 \}$.  The connection $D$ is said to be regular if for all $u\in TM$, we have $ T_u (\T M)=V_u (\T M)\oplus H_u (\T M).$ For a regular connection $D$,  the vector bundle
   maps $  \rho |_{H(\T M)}$ and $K |_{V(\T M)}$
 are isomorphisms. In this case, the map  $\beta:=(\rho |_{H(\T M)})^{-1}$
 is called the horizontal map of $D$. A famous regular connection is the Berwald connection which is defined by \cite{r21}, \cite[Proposition 4.4]{r92}
\begin{equation}\label{Berwald} \gamma {D^\circ}_{h{Z}} \, \overline{W}:={v}[{h}Z,JW], \qquad {D^\circ}_{\gamma\overline{Z}} \, \rho {W}:=\rho[\gamma \overline{Z},{\beta} \overline{W}],
\end{equation}
where $h:= \beta \circ \rho$ is the horizontal projector of $D$ and $v:= I - \beta \circ \rho$ is the  vertical projector of $D$. 

 A nonlinear connection on $M$ is a vector $1$-form $\Gamma$ on $TM$ which is
smooth on $\T M$ and continuous  on $TM$ such that $J \Gamma=J \text{ and } \,\, \Gamma J=-J $ \cite{r21}.  Consequently, the horizontal and vertical projectors
 associated with $\Gamma$ are
given,  respectively, by
 \begin{equation}\label{hor. and ver. proj.}
h:=\frac{1}{2} (I+\Gamma),\quad  \quad v:=\frac{1}{2}
 (I-\Gamma),
\end{equation}
 Moreover,  the curvature of $\Gamma$ is defined by
 $\mathfrak{R}:=-\frac{1}{2}[h,h]$, which can be computed
 using Fr\"{o}licher-Nijenhuis  bracket $[\mathbb{K},\mathbb{L}]$ of two vector 1-forms $\mathbb{K}$ and $\mathbb{L}$ as follows  \cite{r20}:
\begin{eqnarray}\label{brackect of two v.forms}
  [\mathbb{K},\mathbb{L}](W,Z)&=& [\mathbb{K}W,\mathbb{L}Z]+[\mathbb{L}W,\mathbb{K}Z]+\mathbb{K} \mathbb{L}[W,Z]+\mathbb{L} \mathbb{K}[W,Z] \nonumber \\
   && -\mathbb{K}[\mathbb{L} W,Z]-\mathbb{K}[ W,\mathbb{L} Z]-\mathbb{L}[\mathbb{K} W,Z]-\mathbb{L}[ W,\mathbb{K} Z].
\end{eqnarray}
In particular, the vector $2$-form $N_{\mathbb{L}}$ defined by
\begin{equation}\label{Nk}
     N_{\mathbb{L}}(W,Z):= \frac{1}{2}[\mathbb{L},\mathbb{L}](W,Z)=[\mathbb{L}W,\mathbb{L}Z]+\mathbb{L}^{2}[W,Z] -\mathbb{L}[\mathbb{L} W,Z]-\mathbb{L}[ W,\mathbb{L} Z]
\end{equation}
is said to be the Nijenhuis torsion of a vector $1$-form  $\mathbb{L}$. For example, we have $N_{J} = 0$ and $J^{2} =0$ which give
\begin{equation}\label{JJ}
[JW,JZ]=J[W,JZ]+J[JW,Z].
\end{equation}
\begin{defn}\label{regular Finsler}
A Finsler structure \emph{ (}or Finsler metric\emph{) }  on $M$ is a mapping \vspace{-0.3 cm}
$$F: TM \To [0,\infty) $$  such~that{\em:}
 \begin{description}
    \item[(a)] $F $ is  $C^{\infty}$ on  $\T M$ and  $C^{0}$ on $TM$,
    \item[(b)]$F$ is positively homogeneous of degree $1$ in the directional argument $y$, that is
    $\mathcal{L}_{\mathcal{C}} F=F$, where $\mathcal{L}_W$ is the Lie derivative in the direction of $W$,
\item[(c)] the Hilbert $2$-form
    $dd_{J}E$  has a maximal rank, where $E=\frac{1}{2}F^{2} $ is the Finsler energy function.
 \end{description}
The Finsler metric tensor $g$ induced by $\,F\,$ on $\pi^{-1}(TM)$  is defined as follows\emph{ \cite{r94a}}
\begin{equation}\label{g}
g(\rho W,\rho Z):=dd_{J}E(JW,Z), \ \forall W,Z \in
    \mathfrak{X}(TM).
\end{equation}
  In this case, the pair $(M,F)$ is called  a {Finsler manifold} and $F$ is  a regular Finsler metric or simply a Finsler metric. 
   \end{defn}
   According to \cite{Lovas}, the Finsler metric tensor $g$ defined in Definition \ref{regular Finsler} is positive definite.  The main reference for the notion of a conic pseudo-Finsler structure is \cite{Finsler definitions}.
   \begin{defn} A conic sub-bundle of $TM$ is a non-empty open subset $\mathcal{D}  \subset \T M$ that is invariant under scaling of its tangent vectors by positive real numbers and  satisfies $ \pi(\mathcal{D}) = M.$ Assume that, $\forall x \in M, \,\, \mathcal{D}_{x} := D \in T_{x}M$ is a connected   set.  A conic pseudo-Finsler structure $F$ on $M$ is a smooth mapping such that \vspace{-0.3 cm}
$$F: \mathcal{D} \To [0,\infty) ,$$  
$F$ is positively homogeneous of degree $1$ in the directional argument $y$ and the Hilbert $2$-form
    $dd_{J}E$  has a maximal rank. The Finsler metric tensor $g$ induced by $\,F\,$ on $\pi_{\mathcal{D}}^{-1}(TM)$\footnote{$\pi_{\mathcal{D}}$ is the restriction of $\pi$ on  $\mathcal{D}$}  is defined by 
\[g(\rho W,\rho Z):=dd_{J}E(JW,Z), \ \forall \, W,Z \in
    \mathfrak{X}(\mathcal{D}).\]
 The triple $(M,\mathcal{D},F)$ \emph{(}or simply, $(M,F)$\emph{)} is called  a {conic pseudo-Finsler manifold}. 
 \end{defn}
 \begin{rem}
 On a conic domain $\mathcal{D}$, we have again the short exact sequence \emph{(}see,  \emph{\cite{Anastaise})}$$0\longrightarrow
 \pi^{-1}_{\mathcal{D}}(\T M)\stackrel{\gamma}\longrightarrow T{\mathcal{D}} \stackrel{\rho}\longrightarrow
\pi^{-1}_{\mathcal{D}}(\T M)\longrightarrow 0 ,\vspace{-0.1cm}$$ 
with the obvious modifications of the definitions of $\gamma$ and $\rho$.
 \end{rem}
 Further,  the normalized supporting element (or supporting form) is  $\ell :=F^{-1}i_{\overline{\eta}}\:g$ and the angular metric tensor $\hbar := g-\ell \otimes \ell$.
 For a Finsler manifold,  its geodesic spray $G$  satisfies
    $i_{G}\,dd_JE =-d E.$ Moreover, the Barthel connection $\Gamma$ can be written in terms of the geodesic spray as $\Gamma = [J,G]$  \cite{r21}. Another interesting regular connection is the Cartan connection $\nabla$ which is the uniquely determined by  \cite{r92}
   \begin{description}
     \item[(i)]  $ 2g(\nabla _{\gamma\overline{X}}
      \overline{Y}, \overline{Z})
      =\gamma\overline{X}\cdot g( \overline{Y},\overline{Z})+
     g( \overline{Y},\rho [\beta\overline{Z},\gamma\overline{X}])
     +g(\overline{Z},\rho [\gamma\overline{X},\beta \overline{Y}])$,

     \item[(ii)]  $ 2g(\nabla _{\beta\overline{X}}\rho Y,\rho Z)
     = \beta\overline{X}\cdot g( \overline{Y},
                \overline{Z})+
                \beta\overline{Y}\cdot g( \overline{Z},\overline{ X})
                -\beta\overline{Z}\cdot g( \overline{X},
                \overline{Y})$\\
  $ { \qquad \qquad \qquad \quad }
  -g( \overline{X},\rho [\beta\overline{Y},\beta\overline{Z}])+g( \overline{Y},\rho
[\beta\overline{Z},\beta\overline{X}])+g( \overline{Z},\rho
    [\beta\overline{X},\beta\overline{Y}])$.
 \end{description}
\section{ $\phi$-Matsumoto change}
In this section, we give an intrinsic investigation of what we call   generalized  $\phi$-Matsumoto change $F \longrightarrow \widehat{F}$. We find out the relation between the supporting forms ($\ell$ and $\widehat{\ell}$),  the  angular metric tensors ($\hbar$ and $\widehat{\hbar}$),  Finsler metric tensors ($g$ and $\widehat{g}$) and the Cartan torsions ($\mathbf{T}$ and $\widehat{\mathbf{T}}$), corresponding to this change. Moreover, the condition which makes $\widehat{F}$ non-degenerate is derived.
\begin{defn}\emph{\cite{r94a}} \label{concurrent} Let $(M,F)$ be a Finsler manifold.
A non-vanishing $\pi$-vector field  $\overline{\varphi}$ is
called a concurrent $\pi$-vector field if it satisfies
\begin{equation}\label{ch512.eq.1}
     \nabla_{\beta \overline{W}}\,\overline{\varphi}=- \overline{W}={{D}}^{\circ}_{\beta \overline{W}}\,\overline{\varphi} , \quad \quad
     \nabla_{\gamma \overline{W}}\,\overline{\varphi}=0={{D}}^{\circ}_{\gamma \overline{W}}\,\overline{\varphi}.
 \end{equation}
 where $ \nabla \,(\text{respectively,  } {{D}}^{\circ})$ is the Cartan (respectively, Berwald) connection  associated with $F$.
 \end{defn}
Consequently,  the
$\pi$-form $\phi:=i_{\overline{\varphi}}\,g$  associated with $\overline{\varphi}$  has the properties
\begin{equation}
 (\nabla_{\beta \overline{W}}\phi)(\overline{Z})=-g(\overline{W},\overline{Z})=({{D}}^{\circ}_{\beta  \overline{W}}\phi)(\overline{Z}), \quad  (\nabla_{\gamma\overline{W}}\phi)(\overline{Z})=0=({{D}}^{\circ}_{\gamma \overline{W}}\phi)(\overline{Z}).
\end{equation}

\par Let us fix our notation throughout the whole paper:\\
--  $\overline{\varphi}$ denotes a concurrent $\pi$-vector field with respect to $F$,  \\
-- $\phi$ is the $\pi$-form  associated with $\overline{\varphi}$,  \\
--  ${\Phi}:=g(\overline{\varphi},\overline{\eta})=\phi(\overline{\eta})$.
\begin{rem}\label{independent of the directional argument}
For a  Finsler manifold $(M,F)$, a $\pi$-vector field  $\overline{Z} \in \cp$ is independent of the directional argument
 $y$   if,  and only if,  $D^{\circ}_{\gamma \overline{W}}\overline{Z}=0
 $ for all $\overline{W} \in \cp$. Similarly, a scalar (vector) $\pi$-form  $\overline{ \omega }$ is independent of the directional argument
 $y$ if,  and only if, $D^{\circ}_{\gamma \overline{W}}\, \overline{ \omega }=0$ for all $\overline{W} \in \cp$.  Thus, a concurrent $\pi$-vector field $\overline{\varphi}$ and its associated $\pi$-form
$\phi$ are independent of  $y$\emph{ \cite[Theorem 3.7]{r94a}}.  Moreover, we have
\begin{equation}\label{P properties}
 i_{\gamma \overline{\varphi}} \,J  =0, \qquad \qquad d_{J}p^{2}=0, \qquad p^2:=\phi(\overline{\varphi})=g(\overline{\varphi},\overline{\varphi}).
\end{equation}
\end{rem}
\begin{lem}\label{B}  \emph{\cite{r94a, square metric}} Let $(M,F)$ be a Finsler manifold admitting a concurrent $\pi$-vector field $\overline{\varphi}$ with  associated $\pi$-form $\phi$. Then, for all $X \in  \mathfrak{X}(TM)$ and $\overline{W} \in \cp$ we have:
\begin{description}
   \item[(a)]$d_{J}\Phi(\gamma \overline{W})= 0, \quad  d_{J}\Phi(\beta \overline{W})=d \Phi(\gamma \overline{W}) ={{D}}^{\circ}_{\gamma \overline{W}}\Phi = \phi(\overline{W}),$\\
   $ d\Phi(X)=\phi(K X)-F\ell({\rho X})$,
        \item[(b)]$d_{J}\,F(\gamma \overline{W})= 0,  \quad  d_{J}F(\beta \overline{W})=d F(\gamma \overline{W}) = {{D}}^{\circ}_{\gamma \overline{W}}\,F= \ell(\overline{W}),$\\
        $dF(X)=dF(\gamma K X)=\ell(K X),$
             \item[(c)]$d_{h}\,\Phi(\beta \overline{W})=d \Phi(\beta \overline{W})= {{D}}^{\circ}_{\beta \overline{W}}\,\Phi=-F\,\ell(\overline{W}) $, \quad $d \Phi(G)=-F^2$,\\
             $(D^{\circ}_{G} \,\phi)(\overline{W})= -g(\overline{W},\overline{\eta})=-F\,\ell(\overline{W})$,
        \item[(d)]$d_{h}\,F(\beta \overline{W})= d F(\beta \overline{W})={{D}}^{\circ}_{\beta \overline{W}}\,F=0 $,
       \item[(e)] $({{D}}^{\circ}_{\gamma \overline{W}} \,\ell)(\overline{Z})=(\nabla_{\gamma \overline{W}} \,\ell)(\overline{Z})=F^{-1} \hbar(\overline{W},\overline{Z}),\quad \,(D^{\circ}_{G}\, \ell)(\overline{W})= (\nabla_{G}\, \ell)(\overline{W}) =0$,
         \item[(f)]   $dd_{J}E(\gamma \overline{W}, \beta \overline{Z})= g(\overline{W},\overline{Z}), \quad \rho[G,X]=D^{\circ}_G \rho X-K X$,
         \item[(g)] for an arbitrary smooth function $f$ of two variables $ F \text{ and }\Phi$, we have
   \begin{equation}\label{f,b}
   {D^\circ_{\gamma \overline{W}}}f(F,\Phi)=d_Jf(\beta \overline{W})=\frac{\partial f}{\partial F} \, \ell(\overline{W})
        +\frac{\partial f}{\partial \Phi} \, \phi(\overline{W}).
   \end{equation}
     \end{description}
  \end{lem}
\begin{defn}  Let $(M,F)$ be a Finsler manifold admitting a concurrent $\pi$-vector field $\overline{\varphi}$ with the associated $\pi$-form $\phi$. Consider the following change \,
 \begin{equation}\label{change}
\widehat{F}(x,y)=\frac{F^2(x,y)}{F(x,y)-\Phi(x,y)},
\end{equation}
with ${\Phi }=g(\overline{\varphi},\overline{\eta})=\phi(\overline{\eta})$ defined on the conic domain   $\mathcal{D}:=\{(x,y) \in TM \, | \,F(x,y) > \Phi(x,y)\}$.  If $\widehat{F}$ defines a conic pseudo-Finsler structure on $M$,  then $\widehat{F}$ is referred to as   generalized $\phi$-Matsumoto metric.
\end{defn}
 If, in particular,  $F$ is a Riemannian  metric and $\phi$ is an $1$-form on $M$,  then  $\widehat{F}$ is called a Matsumoto metric.  Also,  when $F$ is a Finsler  metric and $\phi$ is an $1$-form on $M$,  $\widehat{F}$ is called a  generalized Matsumoto metric.
The geometric objects corresponding to $\widehat{F}$ will be denoted by hat letters such as $\widehat{\beta},\,\widehat{\ell}, \, \widehat{g},...,$ etc. Thus, clearly we get
\begin{eqnarray}\label{G dot F}
                                                 G \cdot \widehat{F} &=& d\widehat{F}(G)= \frac{F(F-2\Phi)}{(F-\Phi)^2}\,dF(G)+\frac{F^2}{(F-\Phi)^2}\, d\Phi(G)=-\frac{F^4}{(F-\Phi)^2},\\ \label{X dot F}
                                                     W \cdot \widehat{F}&=& d\widehat{F}(W)= \frac{F(F-2\Phi)}{(F-\Phi)^2}\,dF(W)+\frac{F^2}{(F-\Phi)^2}\, d\Phi(W).
 \end{eqnarray}
\begin{prop}\label{hh1}Let $(M,F)$ be a Finsler manifold admitting concurrent $\pi$-vector field $\overline{\varphi}$.  Under the $\phi$-Matsumoto change {\em (\ref{change})}, we have:  \,
\begin{description}\item[(1)] The total derivative of the Finsler energy functions $d\widehat{E}$ and $dE$ are related by
\begin{eqnarray}\label{AB}
  d\widehat{E}
 &=& \frac{F^2(F-2\Phi)}{(F-\Phi)^3}\,dE+\frac{F^4}{(F-\Phi)^3}\, d\Phi .\end{eqnarray}
  \item[(2)] The supporting forms $\widehat{\ell}$ and $\ell$ are related by
  \begin{equation}\label{ell}
     \widehat{\ell} =\frac{F(F-2\Phi)}{(F-\Phi)^2}\,\ell
     +\frac{F^2}{(F-\Phi)^2}\, \phi .
  \end{equation}
  \item[(3)] The vertical counterpart for Berwald connections ${{\widehat{D}}}^{\circ}_{\gamma \overline{W}}\, \overline{Z}$ and ${{D}}^{\circ}_{\gamma \overline{W}}\, \overline{Z}$ are related by
  \begin{equation}\label{invariant v. Berwald connection}
  {{\widehat{D}}}^{\circ}_{\gamma \overline{W}}\, \overline{Z}= {{D}}^{\circ}_{\gamma \overline{W}}\, \overline{Z}.
  \end{equation}
  \item[(4)] The angular metric tensors $\widehat{\hbar}$ and $\hbar$ are related by
  \begin{eqnarray}
    \widehat{\hbar}(\overline{W},\overline{Z}) &=&  \frac{F^2(F-2\Phi)}{(F-\Phi)^3}\hbar(\overline{W},\overline{Z})
    +\frac{2F^4}{(F-\Phi)^4}\phi(\overline{W})\,\phi(\overline{Z})
   +\frac{2 \Phi^2 F^2}{(F-\Phi)^4}\,\ell(\overline{W})\,\ell(\overline{Z}) \nonumber\\
    &&-\frac{2 \Phi F^3}{(F-\Phi)^4}\set{\phi(\overline{W})\,\ell(\overline{Z})+\phi(\overline{Z})\,\ell(\overline{W})}.
      \end{eqnarray}
\end{description}
\end{prop}
\begin{proof}   Under the $\phi$-Matsumoto change {(\ref{change})}: \\
\noindent\textbf{ (1)} We have \begin{eqnarray}
  d\widehat{E}(W)&=&\widehat{F}\,d\widehat{F}(W)=\frac{F^2}{F-\Phi} \set{\frac{F(F-2\Phi)}{(F-\Phi)^2}\,dF(W)+\frac{F^2}{(F-\Phi)^2}\, d\Phi(W)} \nonumber \\
 &=& \frac{F^3(F-2\Phi)}{(F-\Phi)^3}\,dF(W)+\frac{F^4}{(F-\Phi)^3}\, d\Phi(W).\end{eqnarray}
\noindent \textbf{(2)} Taking into account  Lemma \ref{B} \textbf{(a)}, \textbf{(b)}, it follows that
            \begin{eqnarray*}
                \widehat{\ell}(\overline{W})&=& d_J \widehat{F}(\widehat{\beta} \overline{W})=d_J \widehat{F}({\beta} \overline{W})\\
                &=& \frac{\partial \widehat{F}}{\partial F}\, d_J F(\beta \overline{W})
                 +\frac{\partial \widehat{F}}{\partial \Phi} \, d_J \Phi(\beta \overline{W})   \frac{F(F-2\Phi)}{(F-\Phi)^2}\,\ell(\overline{W})
     +\frac{F^2}{(F-\Phi)^2}\, \phi(\overline{W}).
               \end{eqnarray*}
\noindent \textbf{(3)}  Since the difference between the horizontal maps ${\widehat{\beta}}$ and $\beta$ is a vertical vector field,  that is,  $\widehat{\beta}=\beta + \gamma \overline{\mu}$, for some $\overline{\mu}\in \cp$.  Using the facts that $\rho \circ \gamma =0$,   the vertical distribution is completely integrable and
 ${{{D}}}^{\circ}_{\gamma \overline{W}}\overline{Z}=\rho[\gamma
\overline{W}, {\beta} \overline{Z}]$ \cite{TahaSoleiman22}, hence we get
   $$ {{\widehat{D}}}^{\circ}_{\gamma \overline{W}}\overline{Z}=\rho[\gamma
\overline{W}, \widehat{\beta}\, \overline{Z}]=\rho[\gamma \overline{W}, {\beta} \overline{Z}]
+\rho[\gamma \overline{W}, \gamma \overline{\mu}]=\rho[\gamma \overline{W}, {\beta} \overline{Z}]={{{D}}}^{\circ}_{\gamma \overline{W}}\overline{Z}.$$
\noindent \textbf{(4)} Using items \textbf{(2)}, \textbf{(3)}  above, Lemma \ref{B} \textbf{(a)}, \textbf{(b)}, \textbf{(e)}, together with Definition~\ref{concurrent},  one can show that
\begin{eqnarray*}
 \widehat{\hbar}(\overline{W},\overline{Z})&=&\widehat{F}\left({{\widehat{D}}}^{\circ}_{\gamma \overline{W}} \,\widehat{\ell}\,\right)(\overline{Z})=
 \widehat{F}\left({{D}}^{\circ}_{\gamma \overline{W}} \,\widehat{\ell}\,\right)(\overline{Z}) \\
   &=& \widehat{F}\, {{D}}^{\circ}_{\gamma \overline{W}}\left(\frac{F(F-2\Phi)}{(F-\Phi)^2}\,\ell(\overline{Z})
     +\frac{F^2}{(F-\Phi)^2}\, \phi(\overline{Z})\right)\\
   &=& \widehat{F}\, \set{\left({{D}}^{\circ}_{\gamma \overline{W}} \,\frac{F(F-2\Phi)}{(F-\Phi)^2}\right)\,\ell(\overline{Z})+\left({{D}}^{\circ}_{\gamma \overline{W}} \,\frac{F^2}{(F-\Phi)^2}\right)\, \phi(\overline{Z})}\\
   && +\widehat{F}\,\set{\frac{F(F-2\Phi)}{(F-\Phi)^2}\,\left({{D}}^{\circ}_{\gamma \overline{W}} \,\ell\right)(\overline{Z})+\frac{F^2}{(F-\Phi)^2}
   \left({{D}}^{\circ}_{\gamma \overline{W}} \, \phi\right)(\overline{Z})}\\
   &=& \frac{F^2}{F-\Phi}\, \{-\frac{2 \Phi^2}{(\Phi-F)^3}\,\ell(\overline{W})\ell(\overline{Z})+ \frac{2 \Phi F}{(\Phi-F)^3}\, \phi(\overline{W})\,\ell(\overline{Z})\\
   &&  + \frac{2 \Phi F}{(\Phi-F)^3}\,\ell(\overline{W})\phi(\overline{Z}) + \frac{2 F^2}{(F-\Phi)^3}\, \phi(\overline{W}) \, \phi(\overline{Z}) \}\\
   && +\frac{F^2}{F-\Phi}\,\frac{F(F-2\Phi)}{(F-\Phi)^2}\,F^{-1} \,\hbar(\overline{W},\overline{Z}).
\end{eqnarray*}
Hence, the result follows.
\end{proof}
\begin{prop}\label{gg2}The Finsler metric tensor $\widehat{g}$ associated with the generalized $\phi$-Matsumoto metric $\widehat{F}$ is given by  \,
\begin{eqnarray}\label{gg3}
    \widehat{g}(\overline{W},\overline{Z})&=&  \frac{F^2(F-2\Phi)}{(F-\Phi)^3}g(\overline{W},\overline{Z})+\frac{3F^4}{(F-\Phi)^4}\phi(\overline{W})\,\phi(\overline{Z})
    +\frac{F^2\,\Phi(4\Phi-F)}{(F-\Phi)^4}\,\ell(\overline{W})\,\ell(\overline{Z}) \nonumber\\
    &&+\frac{F^3(F-4\Phi)}{(F-\Phi)^4}\set{\phi(\overline{W})\,\ell(\overline{Z})+\phi(\overline{Z})\,\ell(\overline{W})}.
    \end{eqnarray}
\end{prop}
   \begin{proof}
 From the definition of the angular metric tensor $\widehat{\hbar}:=\widehat{g}-\widehat{\ell} \otimes \widehat{\ell}$ and using Proposition \ref{hh1} (2), (4), we obtain
    \begin{eqnarray*}
    \widehat{g}(\overline{W},\overline{Z}) &=&\frac{F^2(F-2\Phi)}{(F-\Phi)^3}\hbar(\overline{W},\overline{Z})
    +\frac{2F^4}{(F-\Phi)^4}\phi(\overline{W})\,\phi(\overline{Z}) \nonumber\\
    &&
   +\frac{2 \Phi^2 F^2}{(F-\Phi)^4}\,\ell(\overline{W})\,\ell(\overline{Z})-\frac{2 \Phi F^3}{(F-\Phi)^4}\set{\phi(\overline{W})\,\ell(\overline{Z})+\phi(\overline{Z})\,\ell(\overline{W})}\\
    &&+ \set{\frac{F(F-2\Phi)}{(F-\Phi)^2}\,\ell(\overline{W})
     +\frac{F^2}{(F-\Phi)^2}\, \phi(\overline{W})}\\
    && \, \times
     \set{\frac{F(F-2\Phi)}{(F-\Phi)^2}\,\ell(\overline{Z})
     +\frac{F^2}{(F-\Phi)^2}\, \phi(\overline{Z})} \\
     &=&   \frac{F^2(F-2\Phi)}{(F-\Phi)^3}g(\overline{W},\overline{Z})+\frac{3F^4}{(F-\Phi)^4}\phi(\overline{W})\,\phi(\overline{Z})
    +\frac{F^2\,\Phi(4\Phi-F)}{(F-\Phi)^4}\,\ell(\overline{W})\,\ell(\overline{Z}) \nonumber\\
    &&+\frac{F^3(F-4\Phi)}{(F-\Phi)^4}\set{\phi(\overline{W})\,\ell(\overline{Z})+\phi(\overline{Z})\,\ell(\overline{W})}.
    \end{eqnarray*}
    \vspace*{-1.1 cm}\[\qedhere\]
      \end{proof}
Consequently, for any $\pi$-vector field $\overline{\mu}$ and  $W \in \mathfrak{X}(TM)$, we get
\begin{eqnarray}\label{AB3}
\widehat{g}(\overline{\mu},\rho W)
 &=& \frac{F^2(F-2\Phi)}{(F-\Phi)^3}g(\overline{\mu},\rho W)+\frac{F^3(F-4\Phi)}{(F-\Phi)^4}\set{\phi(\overline{\mu})\,\ell(\rho W)+\phi(\rho W)\,\ell(\overline{\mu})} \nonumber\\ &&
 +\frac{3F^4}{(F-\Phi)^4}\phi(\overline{\mu})\,\phi(\rho W)
    +\frac{F^2\,\Phi(4\Phi-F)}{(F-\Phi)^4}\,\ell(\overline{\mu})\,\ell(\rho W).
\end{eqnarray}
\begin{thm}
Let $(M, F)$ be a Finsler manifold admitting a concurrent $\pi$-vector field $\overline{\varphi}$ with associated  $\pi$-form $\phi$. The function $\widehat{F}$ defined by \eqref{change}  is  a conic pesudo-Finsler structure if and only if \begin{equation}
\label{Eq:g_tilde_non_degenerate}
F(1+2p^2)-3\Phi\neq0.
\end{equation}
That is,  the Finsler metric tensor $\widehat{g}$ of $\widehat{F}$ is non-degenerate if and only if the function $F(1+2p^2)-3\Phi$  does not vanish.
\end{thm}
 \begin{proof}
 The metric $\widehat{g}$ associated with $\widehat{F}$ is non-degenerate if and only if    $$\widehat{g}(\overline{W}, \overline{Z})=0\,\, \, \,\forall\, \overline{W}\in\cp \implies \overline{Z}=0.$$
 Assume that  $\widehat{g}(\overline{W}, \overline{Z})=0,\,\, \, \forall\, \overline{W}\in\cp$.  Then, relation \eqref{gg3} gives rise to
\begin{eqnarray}\label{sunsb}
   0 &=&  \frac{F^2(F-2\Phi)}{(F-\Phi)^3}g(\overline{W},\overline{Z})+\frac{3F^4}{(F-\Phi)^4}\phi(\overline{W})\,\phi(\overline{Z})
    +\frac{F^2\,\Phi(4\Phi-F)}{(F-\Phi)^4}\,\ell(\overline{W})\,\ell(\overline{Z}) \nonumber\\
    &&+\frac{F^3(F-4\Phi)}{(F-\Phi)^4}\set{\phi(\overline{W})\,\ell(\overline{Z})+\phi(\overline{Z})\,\ell(\overline{W})}.
    \end{eqnarray}
 Setting $\overline{W}=\overline{\varphi}$ in \eqref{sunsb}, noting that $\ell(\overline{\varphi})=\frac{\Phi}{F}$ and
$\phi(\overline{\varphi})=p^2$, one can show that
\begin{equation}\label{12q}
  \chi_1 \,\ell(\overline{Z}) + \Upsilon_1\, \phi({\overline{Z}})=0,
\end{equation}
where
$$  \chi_1 := \frac{F (4 \Phi-F) \left(\Phi^2-F^2 p^2\right)}{(\Phi-F)^4}, \quad
  \Upsilon_1 :=\frac{F^2 \left(-2 \Phi^2-2 \Phi F+F^2 \left(3 p^2+1\right)\right)}{(\Phi-F)^4}.$$
Similarly,  setting by $\overline{W}=\overline{\eta}$ in \eqref{sunsb}, taking into account $\ell(\overline{\eta})=F$ and
$\phi(\overline{\eta})=\Phi$, we obtain
\begin{equation}\label{13q}
  \chi_2\, \ell(\overline{Z}) + \Upsilon_2\, \phi({\overline{Z}})=0,
\end{equation}
with
 $$ \chi_2 := \frac{F^3 (F-2 \Phi)}{(F-\Phi)^3},\quad
  \Upsilon_2 :=\frac{F^4}{(F-\Phi)^3}.$$
The system of the  algebraic equations (\ref{12q}) and (\ref{13q}), in $ \ell \text{ and } \phi$, has a non-trivial solution if and only if
$$\frac{F^6\, ( F \left(2 p^2+1\right)-3 \Phi)}{(\Phi-F)^6}=0. $$
Hence, as $\widehat{F}(x,y)= \frac{F^2(x,y)}{F(x,y)-\Phi(x,y)}\neq0$ over $\T M$, then we conclude that $F \left(2 p^2+1\right)-3 \Phi=0$.\\
Consequently,
$$\overline{Z}\neq0 \quad \Longleftrightarrow \quad F \left(2 p^2+1\right)-3 \Phi=0.$$
Therefore,  $\overline{Z}=0$ if and only if  the Finsler structure $F$ and the $\pi$-form $\Phi$ satisfy the condition
$$F \left(2 p^2+1\right)-3 \Phi\neq0.$$
This means that
the metric tensor $\widehat{g}$ of $\widehat{F}$ is non-degenerate if and only if the condition \eqref{Eq:g_tilde_non_degenerate} is satisfied. Hence, the proof is complete.
\end{proof}
     Form now on, we consider that the  generalized  $\phi$-Matsumoto metric  $\widehat{F}$,  defined by \eqref{change},   satisfies the condition \eqref{Eq:g_tilde_non_degenerate}.
\begin{prop} Let $(M,F)$ be a Finsler manifold admitting a concurrent $\pi$-vector field $\overline{\varphi}$.   Under the $\phi$-Matsumoto change {\em (\ref{change})},  the Cartan torsion tensor $\widehat{\mathbf{T}}$ of  $\widehat{F}$ can be written in terms of the Cartan torsion tensor ${\mathbf{T}}$ of $F$ as follows
\begin{eqnarray*}\label{tt3}
    2\widehat{\mathbf{T}}(\overline{W},\overline{Y},\overline{Z})&=&  \frac{2F^2(F-2\Phi)}{(F-\Phi)^3}\mathbf{T}(\overline{W},\overline{Y},\overline{Z})
    \\
    &&+\frac{F\,\Phi(4\Phi-F)}{(F-\Phi)^4}\,\set{\hbar(\overline{W},\overline{Z})\,\ell(\overline{Y})
    +\hbar(\overline{Y},\overline{Z})\,\ell(\overline{W})} \nonumber\\
    &&+\frac{F^2(F-4\Phi)}{(F-\Phi)^4}\set{\hbar(\overline{W},\overline{Z})\,\phi(\overline{Y})
    +\hbar(\overline{Y},\overline{Z})\,\phi(\overline{W})} \\
    &&+\left({D^\circ_{\gamma \overline{Z}}}\, \frac{F^2(F-2\Phi)}{(F-\Phi)^3}\right)g(\overline{W},\overline{Y})+{D^\circ_{\gamma \overline{Z}}}\left(\frac{3F^4}{(F-\Phi)^4}\,\phi(\overline{W})\,\phi(\overline{Y})\right)\\
    &&
    +\left({D^\circ_{\gamma \overline{Z}}}\frac{F^2\,\Phi(4\Phi-F)}{(F-\Phi)^4}\right)\,\,\ell(\overline{W})\,\ell(\overline{Y}) \nonumber\\
    &&+\left({D^\circ_{\gamma \overline{Z}}}\frac{F^3(F-4\Phi)}{(F-\Phi)^4}\right)\,\set{\phi(\overline{W})\,\ell(\overline{Y})+\phi(\overline{Y})\,\ell(\overline{W})}.
        \end{eqnarray*}
\end{prop}
\begin{proof}
Using the expression of the metric $\widehat{g}$ obtained in Proposition \ref{gg2}, taking into account   the fact that
    $({D^\circ_{\gamma \overline{Z}}}g)(\overline{W},\overline{Y})=2\mathbf{T}(\overline{W},\overline{Y},\overline{Z})$,  the result follows.
\end{proof}
\section{The change of the spray and Berwald connection}
Here, we find out the expression of  the geodesic spray $\widehat{G}$ of the generalized  $\phi$-Matsumoto metric $\widehat{F}$ in terms of the geodesic spray $G$ of $F$.  Consequently,  we prove that the two geodesic sprays $G$ and $\widehat{G}$ can not be projectively related.
Moreover,  the relationship between the two Barthel connections $\Gamma$ and $\widehat{\Gamma}$ is obtained,  as well as the relations between Barthel curvature tensors $\Re $ and $\widehat{\Re}$  and Berwald connections $D^\circ$ and  $\widehat{D^\circ}$ are derived.  We conclude that the $\pi$-vector field $\overline{\varphi}$ is never concurrent with respect to $\widehat{F}$,  then we end this section by providing a condition  that makes   $\overline{\varphi}$ concurrent with respect to  $\widehat{F}$.
\begin{thm}\label{th.22} Let $(M,F)$ be a Finsler manifold admitting a concurrent $\pi$-vector field $\overline{\varphi}$ with  associated $\pi$-form $\phi$.  If $G$ is the geodesic spray of $F$,  then the geodesic spray $\widehat{G}$ of the generalized $\phi$-Matsumoto metric $\widehat{F}$  is given by
$$\widehat{G}=G-f_1\,{\mathcal{C}}+f_2\,\gamma \overline{\varphi},$$
where $f_1:=\frac {F(4\Phi-F)} {F(1+2p^2)-3\Phi} \quad \text{and} \quad f_2:=\frac{2F^3}{F(1+2p^2)-3\Phi}$.
\end{thm}
\begin{proof}
Since the geodesic spray $\widehat{G}$  of $\widehat{F}$ satisfies \cite{r21} $$-d\widehat{E}= \frac{1}{2}\,i_{\widehat{G}} dd_{J} \widehat{F}^2,$$where $d\widehat{E}$ is expressed in \eqref{AB} in terms of $dE$.
 Due  to the fact that the difference between two sprays is a vertical vector field (i.e. $\widehat{G}=G+\gamma \overline{\mu}$, for some
$\pi$-vector field $\overline{\mu}$), we get
\begin{equation}\label{ch5eq}
-d\widehat{E}(W)=\frac{1}{2}\,i_{G+\gamma
\overline{\mu}}\,dd_{J} \widehat{F}^2(W)=\frac{1}{2}i_{G}\,dd_{J}\widehat{F}^2(W)+\frac{1}{2}i_{\gamma \overline{\mu}}\,dd_{J}\widehat{F}^2(W),
 \end{equation}
 Using  $\beta \overline{\eta}=G$ and $W=hW+vW=\beta\rho W+\gamma K W$, together with Lemma \ref{B}, we obtain
 \begin{eqnarray}\label{AB1}
   \frac{1}{2}\,i_{G}\,dd_{J}\widehat{F}^{2}(W)&=&\frac{1}{2}\{dd_{J}\widehat{F}^{2}(\beta \overline{\eta} ,W)\}=\frac 1 2 \set{G \cdot d_J\widehat{F}^{2}(W) - W \cdot d_J\widehat{F}^{2}(G)-d_J\widehat{F}^{2}[G,W]} \nonumber \\
 &=& \frac 1 2 \set{G \cdot (2\widehat{F} \widehat{\ell}(\rho W))-W \cdot(2 \widehat{F} \widehat{\ell}(\overline{\eta}))
 -2\widehat{F}\widehat{\ell}(\rho[G,W])} \nonumber \\
 &=&{(G \cdot \widehat{F})\,\widehat{\ell}(\rho W)+\widehat{F}\, G \cdot \widehat{\ell}(\rho W)-(W \cdot \widehat{F}^2)
 -\widehat{F}\,\widehat{\ell}(\rho[G,W])}.
\end{eqnarray}
Taking into account Lemmas \ref{B}, Proposition \ref{hh1} (2), Formula \eqref{G dot F} and \eqref{X dot F}, relation \eqref{AB1} reduces to
\begin{eqnarray*}
             \frac{1}{2}\,i_{G}\,dd_{J}\widehat{F}^{2}(W)
 &=&-\frac{F^4}{(F-\Phi)^2} \left(\frac{F(F-2\Phi)}{(F-\Phi)^2}\, \ell(\rho {W})+ \frac{F^2}{(F-\Phi)^2}\, \phi(\rho {W})\right) \nonumber\\
 &&+\frac{F^2}{(F-\Phi)}\, G \cdot\left( \frac{F(F-2\Phi)}{(F-\Phi)^2}\, \ell(\rho {W})+  \frac{F^2}{(F-\Phi)^2}\, \phi(\rho {W})\right) \nonumber\\
 && -\frac{2F^2}{(F-\Phi)}\,\left(\frac{F(F-2\Phi)}{(F-\Phi)^2} dF(W)+\frac{F^2}{(F-\Phi)^2} d\Phi(W)\right) \nonumber\\
  &&-\frac{F^2}{(F-\Phi)}\,\, \left(\frac{F(F-2\Phi)}{(F-\Phi)^2}\,\, \ell(\rho[G,W])+ \frac{F^2}{(F-\Phi)^2}\, \phi(\rho[G,W])\right) \nonumber\\
  \end{eqnarray*}
  \begin{eqnarray}\label{AB2}
 &=& - \frac{F^5 (F-4 \Phi)}{(\Phi-F)^4} \ell(\rho {W})-
  \frac{3 F^6}{(\Phi-F)^4}\, \phi(\rho {W}) \nonumber \nonumber\\
 &&-\frac{F^3(F-2\Phi)}{(F-\Phi)^3}\,dF(W)-\frac{F^4}{(F-\Phi)^3}\, d\Phi(W).
\end{eqnarray}
Plugging  the relations \eqref{AB2} and \eqref{AB3} (which expresses $\frac{1}{2}\,i_{\gamma \overline{\mu}}\,dd_{J}\widehat{F}^{2}(W)$ as $\widehat{g}(\overline{\mu},\rho W)$) into Equation (\ref{ch5eq}), after some calculation, it follows that \,
\begin{eqnarray*}
            &&- \frac{F^3(F-2\Phi)}{(F-\Phi)^3}\,dF(W)-\frac{F^4}{(F-\Phi)^3}\, d\Phi(W)\\
             &=&- \frac{F^5 (F-4 \Phi)}{(\Phi-F)^4} \ell(\rho {W})-
  \frac{3 F^6}{(\Phi-F)^4}\, \phi(\rho {W})-\frac{F^3(F-2\Phi)}{(F-\Phi)^3}\,dF(W)-\frac{F^4}{(F-\Phi)^3}\, d\Phi(W)\\
 &&+ \frac{F^2(F-2\Phi)}{(F-\Phi)^3}g(\overline{\mu},\rho W)+\frac{3F^4}{(F-\Phi)^4}\phi(\overline{\mu})\,\phi(\rho W)
    +\frac{F^2\,\Phi(4\Phi-F)}{(F-\Phi)^4}\,\ell(\overline{\mu})\,\ell(\rho W) \nonumber\\
    &&+\frac{F^3(F-4\Phi)}{(F-\Phi)^4}\set{\phi(\overline{\mu})\,\ell(\rho W)+\phi(\rho W)\,\ell(\overline{\mu})}.
\end{eqnarray*}
Using the non-degeneracy property of Finsler metric $g$, the above relation reduces to
\begin{eqnarray}\label{ch52.eq.5}
   \frac{F^2(F-2\Phi)}{(F-\Phi)^3}\,\overline{\mu}&=&
   \set{\frac{F^4(F-4\Phi)}{(F-\Phi)^4} -\frac{F^2(F-4\Phi)}{(F-\Phi)^4} \phi({\overline{\mu}})-\frac{F\,\Phi(4\Phi-F)}{(F-\Phi)^4}\ell(\overline{\mu})}\overline{\eta} \nonumber \\
   && +   \set{\frac{3F^6}{(F-\Phi)^4}-\frac{3F^4}{(F-\Phi)^4}\,\phi(\overline{\mu})
   -\frac{F^3(F-4\Phi)}{(F-\Phi)^4}\,\ell(\overline{\mu})}\overline{\varphi},\,
\end{eqnarray}
where $\ell(\overline{\mu})$ and $\phi(\overline{\mu})$ are  determined by the following two equations
\begin{eqnarray}\label{mm}
    &&{F^2(F-2\Phi)}\ell(\overline{\mu})+{F^3}\phi(\overline{\mu})
    ={F^5} ,\\
        &&{F^4(-4\Phi^2+\Phi\,F+3F^2\,p^2)} =
    \nonumber\\
    &&{F\,(4\Phi-F)(\Phi^2-F^2\,p^2)}\ell(\overline{\mu})
    +{F^2(2\Phi^2-2\Phi\,F+F^2(1+3p^2))}\phi(\overline{\mu}).
\end{eqnarray}
Thus, the condition \eqref{Eq:g_tilde_non_degenerate} leads to
\begin{equation}
  \ell(\overline{\mu}) = \frac{F^2 (2 \Phi-F)}{3 \Phi-F \left(2 p^2+1\right)}, \quad
  \phi(\overline{\mu}) =\frac{F \left(-4 \Phi^2+\Phi F+2 F^2 p^2\right)}{-3 \Phi+2 F p^2+F}.
\end{equation}
Consequently, in view of Equation (\ref{ch52.eq.5}) and the assumption $\widehat{G}=G+\gamma \overline{\mu}$, it follows that
the geodesic sprays  $\widehat{G}$ is given by
\begin{eqnarray}\label{geodesic sprays change}
\widehat{G}=G-\frac {F(4\Phi-F)} {F(1+2p^2)-3\Phi} \,{\mathcal{C}}+\frac{2F^3}{F(1+2p^2)-3\Phi}\,\gamma \overline{\varphi}.
\end{eqnarray}
Hence, the proof is complete.
\end{proof}
\begin{thm}
Let $(M,F)$ be a Finsler manifold admitting concurrent $\pi$-vector field $\overline{\varphi}$.   Under the $\phi$-Matsumoto change {\em (\ref{change})},  the geodesic sprays $G$ and  $\widehat{G}$ can never be projectively related.
\end{thm}
\begin{proof} Let $G$ and  $\widehat{G}$ be the geodesic sprays of $F$ and the generalized $\phi$-Matsumoto metric $ \widehat{F}$, respectively.
The Finsler metrics $F$ and $ \widehat{F}$ are projectively related means that $$\widehat{G} = G - 2 A(x,y) \,\mathcal{C},$$ with a projective factor $A(x,y) $ being positively homogeneous function of degree $1$ in the directional argument $y$.  As $F $ is a non-zero function, in view of the relation \eqref{geodesic sprays change}, we get $G$ and  $\widehat{G}$ are projectively related if and only if $\gamma \overline{\varphi}=0.$ Thus, $A(x,y) =\frac{1}{2} f_{1}$, which leads to $\overline{\varphi}=0$.  This contradicts our assumption that the $\pi$-vector field $\overline{\varphi}$ is everywhere nonzero.
\end{proof}
\begin{prop}\label{th.barthel} The Barthel connection $\widehat{\Gamma}$ associated with the generalized $\phi$-Matsumoto metric $\widehat{F}$ can be expressed in terms of the Barthel connection $\Gamma$  associated with $F$ as follows
\begin{equation*}\label{ch52.eq.4}
         \widehat{\Gamma} =\Gamma -f_1\,J-d_J f_1 \otimes \gamma \overline{\eta} + d_J f_2 \otimes \gamma \overline{\varphi}.
\end{equation*}
  Consequently, the horizontal map $\widehat{\beta}$  associated with  the  $\widehat{F}$ has the form
$$\widehat{\beta}\overline{W} ={\beta}\overline{W} -{\frac 1 2}\set{f_1\,\gamma \overline{W}+d_J f_1(\beta\overline{W}) \, \gamma \overline{\eta} - d_J f_2 (\beta\overline{W})\, \gamma \overline{\varphi}}.$$
 \end{prop}
\begin{proof}
In view of Theorem \ref{th.22} and \eqref{eq. 1} along with \eqref{J properties}, we obtain \,
\begin{eqnarray*}
  \widehat{\Gamma} &=& [J,\widehat{G}]
  = \left[J, G-f_1\, \gamma \overline{\eta}+f_2\, \gamma \overline{\varphi} \right]
   =[J,G]+[f_1\, \gamma \overline{\eta}-f_2\, \gamma \overline{\varphi},J]\\
  &=&  [J,G]+f_1[\gamma \overline{\eta},J]+df_1\wedge i_{\gamma\overline{\eta}}\,J-d_{J}f_1  \otimes \gamma \overline{\eta} -f_2[\gamma \overline{\varphi},J]-df_2\wedge i_{\gamma\overline{\varphi}}\,J+d_{J}f_2 \otimes \gamma \overline{\varphi}.
\end{eqnarray*}
From \eqref{J properties} and \eqref{P properties} along with
 \begin{eqnarray*}
 [ \gamma \overline{\varphi}, J]W &=&[\gamma \overline{\varphi}, JW]-J[\gamma \overline{\varphi}, W] \\
  &=&\gamma \{\nabla_{\gamma \overline{\varphi}}\,\rho W-\nabla_{J W}\,\overline{\varphi}\}
  -\gamma\{\nabla_{\gamma \overline{\varphi}}\,\rho W-T(\overline{\varphi},\rho W)\}=0,
  \end{eqnarray*}
 we get
$$\widehat{\Gamma} =\Gamma -f_1\,J-d_J f_1 \otimes \gamma \overline{\eta} + d_J f_2 \otimes \gamma \overline{\varphi}.$$
 Consequently, using the fact that $\Gamma=2\beta \circ \rho -I$, the horizontal map $\widehat{\beta}$  associated with  the generalized $\phi$-Matsumoto metric has the form
$$\widehat{\beta}\,\overline{W} ={\beta}\overline{W} -{\frac 1 2}\set{f_1\,\gamma \overline{W}+d_J f_1(\beta\overline{W}) \, \gamma \overline{\eta} - d_J f_2 (\beta\overline{W})\, \gamma \overline{\varphi}}.$$
\vspace*{-1.1 cm}\[\qedhere\]
\end{proof}
\begin{cor}\label{Coro. hor. and ver. proj.}
The horizontal projector $\widehat{h}$ and vertical projector $\widehat{v}$  associated with  the generalized $\phi$-Matsumoto metric  $\widehat{F}$  can be written,  in terms of the horizontal projector $h$ and vertical projector $v$  associated with $F$,  in the form
\begin{equation}\label{h and v}
\widehat{h}=h+\mathbb{S}, \quad  \widehat{v}=v-\mathbb{S},
\end{equation}
where $\mathbb{S}$ is a semi-basic vector $1$-form \emph{\cite{r21} }given by
\begin{equation}\label{LLL}
  \mathbb{S}:=-\frac 1 2 \set{f_1\,J+d_J f_1 \otimes \gamma \overline{\eta} - d_J f_2 \otimes \gamma \overline{\varphi}}.
\end{equation}
Moreover, the Barthel curvature tensor  $\widehat{\Re}$ associated with the generalized $\phi$-Matsumoto metric {\em (\ref{change})} is determined by
\begin{equation*}
         \widehat{\Re} =\Re -[h,\mathbb{S}]-N_\mathbb{S}.
\end{equation*}
 \end{cor}
\begin{proof} The relation \eqref{h and v} follows from Theorem \ref{th.barthel} and formula~\eqref{hor. and ver. proj.}.
 Now, since the Barthel curvature tensor of $\widehat{F}$ is defined by   $\widehat{\Re}=-\frac 1 2 [\widehat{h},\widehat{h}]$,  the proof follows from  \eqref{h and v} and formula \eqref{Nk} together with the properties of the Fr\"{o}licher-Nijenhuis  bracket. More precisely,   $$\widehat{\Re}=-\frac 1 2\, [h+\mathbb{S},h+\mathbb{S}] = -\frac 1 2 \left([h,h] + [h,\mathbb{S}] +[\mathbb{S},h] + [\mathbb{S},\mathbb{S}]\right) = \Re -[h,\mathbb{S}]-N_\mathbb{S}.$$
 \vspace*{-1.15 cm}\[\qedhere\]
\end{proof}
\begin{prop}\label{th.Berwald} For the generalized $\phi$-Matsumoto metric {\em (\ref{change})}, we have:
\begin{description}
  \item[(1)] The vertical counterpart of Berwald connection can be expressed as
   $$\widehat{D^\circ}_{\gamma\overline{W}} \, \overline{Z}={D^\circ}_{\gamma\overline{W}} \, \overline{Z}.$$
  \item[(2)] The horizontal counterpart of  Berwald connection can be expressed as
  \begin{eqnarray*}
 \widehat{D^\circ}_{\widehat{\beta} \overline{W}} \, {\overline{Z}}&=&   {D^\circ}_{{\beta} \overline{W}} \overline{Z}-
 \frac 1 2\{f_1\,D^\circ_{\gamma \overline{W}}\,\overline{Z}+d_J f_1 ({\beta} \overline{W})\,
    D^\circ_{\gamma \overline{\eta}}\, \overline{Z}\\
    &&-d_J f_1 ({\beta} \overline{W})\, \overline{Z}- d_J f_1 (\beta \overline{Z})\,\overline{W}
    - d_J f_2 ({\beta} \overline{W})\,D^\circ_{\gamma \overline{\varphi}}\, \overline{Z} \} \\
    && +\frac 1 2\set{dd_J f_1 (\gamma \overline{Z},{\beta} \overline{W})\,  \overline{\eta} - dd_J f_2 (\gamma \overline{Z}, {\beta} \overline{W}) \, \overline{\varphi}}.
\end{eqnarray*}
  \end{description}
 \end{prop}
\begin{proof}  Under the  $\phi$-Matsumoto change \eqref{change},   we have:
\begin{description}
  \item[(1)]  The vertical counterpart for Berwald connection ${{D}}^{\circ}_{\gamma \overline{W}} \overline{Z}$ is invariant by \eqref{invariant v. Berwald connection}.
 \item[(2)] Using the facts that $v:=\gamma \circ K$, $h:=\beta \circ \rho$,  the Berwald v-curvature $\widehat{S}^{\circ}=0$, together with formulae \eqref{Berwald}, \eqref{brackect of two v.forms} and \eqref{JJ}, we obtain
 \begin{eqnarray*}
   \gamma \widehat{D^\circ}_{h {W}} \, \rho {Z}&=&\widehat{v}[ \widehat{h}W,JZ]\overset{\eqref{h and v}}{=}(v-\mathbb{S})[(h+\mathbb{S})W,JZ] \\
    &=& v[hW,JZ]+v[\mathbb{S}W,JZ]-\mathbb{S}[hW,JZ]-\mathbb{S}[\mathbb{S}W,JZ] \\
    &\overset{\eqref{LLL}}{=}&  \gamma {D^\circ}_{h {W}} \overline{Z}\\
    &&  -\frac \gamma 2\set{f_1\,K[\,JW,JZ]+d_J f_1 (W)\, K[\, \gamma \overline{\eta},JZ] - d_J f_2 (W)\, K[\,\gamma \overline{\varphi},JZ]}\\
    &&+\frac \gamma 2\set{(JZ \cdot f_1)\,\rho W+(JZ \cdot d_J f_1 (W))\,  \overline{\eta} - (JZ \cdot d_J f_2 (W)) \, \overline{\varphi}}\\
    &&+\frac {\gamma} 2 \set{f_1\,\rho([hW,JZ])+d_J f_1 ([hW,JZ])  \, \overline{\eta} - d_J f_2 ([hW,JZ])  \, \overline{\varphi}}\\
    &=&  \gamma {D^\circ}_{h {W}} \rho {Z}-\frac \gamma 2\{f_1\,D^\circ_{J W}\,\rho Z+d_J f_1 (W)\,
    D^\circ_{\gamma \overline{\eta}}\,\rho Z\\
    &&-d_J f_1 (W)\,\rho Z- d_J f_1 (Z)\,\rho W
    - d_J f_2 (W)\,D^\circ_{\gamma \overline{\varphi}}\, \rho Z \} \\
    && +\frac \gamma 2\set{dd_J f_1 (JZ,W)\,  \overline{\eta} - dd_J f_2 (JZ, W) \, \overline{\varphi}}.
    \end{eqnarray*}
Consequently,
\begin{eqnarray*}
 \widehat{D^\circ}_{\widehat{\beta} \,\overline{W}} \, {\overline{Z}}&=&   {D^\circ}_{{\beta} \overline{W}} \overline{Z}-
 \frac 1 2\{f_1\,D^\circ_{\gamma \overline{W}}\,\overline{Z}+d_J f_1 ({\beta} \overline{W})\,
    D^\circ_{\gamma \overline{\eta}}\, \overline{Z}\\
    &&-d_J f_1 ({\beta} \overline{W})\, \overline{Z}- d_J f_1 (\beta \overline{Z})\,\overline{W}
    - d_J f_2 ({\beta} \overline{W})\,D^\circ_{\gamma \overline{\varphi}}\, \overline{Z} \} \\
    && +\frac 1 2\set{dd_J f_1 (\gamma \overline{Z},{\beta} \overline{W})\,  \overline{\eta} - dd_J f_2 (\gamma \overline{Z}, {\beta} \overline{W}) \, \overline{\varphi}}.
\end{eqnarray*}
\end{description}
\vspace*{-1.1 cm}\[\qedhere\]
\end{proof}
\begin{cor}
 Let $(M,F)$ be a Finsler manifold admitting a concurrent $\pi$-vector field $\overline{\varphi}$.  Under the  $\phi$-Matsumoto change \eqref{change},  the $\pi$-vector field $\overline{\varphi}$ can not be  concurrent with respect to the Finsler metric $\widehat{F}$.
\end{cor}
\begin{proof}
It follows directly from Definition \ref{concurrent}  and Proposition \ref{th.Berwald} as $$\widehat{D^\circ}_{\widehat{\beta} \,\overline{W}} \, {\overline{Z}} \neq {D^\circ}_{{\beta} \overline{W}} \overline{Z}.$$
\vspace*{-1.5 cm}\[\qedhere\]
\end{proof}
Nevertheless, the following result gives a condition  under which the property of the $\pi$-vector field $\overline{\varphi}$ being concurrent is preserved under  $\phi$-Matsumoto change.
\begin{thm}
 Let  $(M,F)$ be a Finsler manifold admitting a concurrent $\pi$-vector field  $\,\overline{\varphi}$.  Let $\widehat{F}$ be the  generalized $\phi$-Matsumoto metric defined by \eqref{change}.  A necessary and sufficient condition for $\overline{\varphi}$ to be concurrent with respect to  $\widehat{F}$ is that
\begin{equation}\label{preserved concurrent}
[d_J f_1 (\beta \overline{X}) - dd_J f_2 (\gamma \overline{\varphi}, {\beta} \overline{X})  ] \overline{\varphi}
- d_J f_1 (\beta \overline{\varphi})  \overline{X}
+dd_J f_1 (\gamma \overline{\varphi},{\beta} \overline{X})\,  \overline{\eta}
=0
\end{equation}
\end{thm}
\begin{proof}
In view of Proposition \ref{th.Berwald}, we have
\begin{equation}\label{h-concurrent}
\widehat{D^\circ}_{\gamma\overline{X}} \, \overline{\varphi}={D^\circ}_{\gamma\overline{X}} \, \overline{\varphi}
\end{equation}
and
\begin{eqnarray}\label{v-concurrent}
 \widehat{D^\circ}_{\widehat{\beta} \,\overline{X}} \, {\overline{\varphi}}&=&   {D^\circ}_{\beta \overline{X}} \, {\overline{\varphi}} + \frac 1 2\set{ f_1\,D^\circ_{\gamma \overline{X}}\,\overline{\varphi}+d_J f_1 ({\beta} \overline{X})\,
    D^\circ_{\gamma \overline{\eta}}\, \overline{\varphi} - d_J f_2 ({\beta} \overline{X})\,D^\circ_{\gamma \overline{\varphi}}\, \overline{\varphi} }\\ \nonumber
 && +\frac 1 2\set{d_J f_1 (\beta \overline{X})   \overline{\varphi}- d_J f_1 (\beta \overline{\varphi})  \overline{X} +dd_J f_1 (\gamma \overline{\varphi},{\beta} \overline{X})\,  \overline{\eta} - dd_J f_2 (\gamma \overline{\varphi}, {\beta} \overline{X}) \, \overline{\varphi}}.
\end{eqnarray}
Now, assume that the condition \eqref{preserved concurrent} is satisfied, then \eqref{v-concurrent} reduces to
\begin{eqnarray}\label{v-concurrent-cond}
\widehat{D^\circ}_{\widehat{\beta}\, \overline{X}} \, {\overline{\varphi}}&=&   {D^\circ}_{\beta \overline{X}} \, {\overline{\varphi}} + \frac 1 2\set{ f_1\,D^\circ_{\gamma \overline{X}}\,\overline{\varphi}+d_J f_1 ({\beta} \overline{X})\,
    D^\circ_{\gamma \overline{\eta}}\, \overline{\varphi} - d_J f_2 ({\beta} \overline{X})\,D^\circ_{\gamma \overline{\varphi}}\, \overline{\varphi}}.
\end{eqnarray}
Since  $\overline{\varphi}$ is a concurrent $\pi$-vector field with respect to  $F$, i.e.,  ${D^\circ}_{{\beta} \overline{X}} \overline{\varphi} =- \overline{X}$ and ${D^\circ}_{\gamma \overline{X}}\,\overline{\varphi}=0$, then, by \eqref{v-concurrent-cond}, $\widehat{D^\circ}_{\widehat{\beta}\, \overline{X}} \, {\overline{\varphi}}=-\overline{X}$. This, together with \eqref{h-concurrent}, imply that $\overline{\varphi}$ is a concurrent $\pi$-vector field with respect to $\widehat{F}$.\\
\par Now, if $\overline{\varphi}$ is concurrent $\pi$-vector field with respect to $\widehat{F}$, i.e., $\widehat{D^\circ}_{\gamma\overline{X}} \, \overline{\varphi} = 0$ and $\widehat{D^\circ}_{\widehat{\beta}\, \overline{X}} \, {\overline{\varphi}}=-\overline{X}$, then, by  \eqref{h-concurrent} and \eqref{v-concurrent-cond}, $$ f_1\,D^\circ_{\gamma \overline{X}}\,\overline{\varphi}+d_J f_1 ({\beta} \overline{X})\,
    D^\circ_{\gamma \overline{\eta}}\, \overline{\varphi} - d_J f_2 ({\beta} \overline{X})\,D^\circ_{\gamma \overline{\varphi}}\, \overline{\varphi}=0.$$
    \vspace*{-1.1 cm}\[\qedhere\]
\end{proof}
\section{Preservation of (almost) rationality }
Now, we give an example of a conic non-Riemannian Finsler metric that admits a concurrent $\pi$-vector field and find $\widehat{F}$.
\begin{example}\label{example}
\em{
Let $M=\Real^3$ and $F$ be a conic Finsler metric  defined by
\begin{equation}\label{Lexample}
F(x,y):=F(x_1 ,x_2, x_3; y_1, y_2, y_3)=\sqrt{x_3^2\left(\frac{x_1^2\,  y_2^2+2y_1y_2}{y_1}\right)^2+y_3^2}
\end{equation}
on the domain $D_{F}=\{(x_1 ,x_2, x_3; y_1, y_2, y_3) \in T \Real^3 \, | \, x_1\neq0, x_3 \neq 0,  y_1 \neq 0, y_2\neq0 \}$.

By straightforward calculations or using the Maple Finsler package \cite{NF_Package},  the non-vanishing components of the Finsler metric tensor $g_{ij}$  are given by
$$\qquad \qquad  \qquad g_{11}=\frac{x_3^2\,   x_1^2\, y_2^3(3x_1^2 y_2+4y_1)}{y_1^4},\quad\quad \qquad g_{12}=-\frac{2 x_3^2\,x_1^2 y_2^2 (2x_1^2\, y_2 +3\,y_1 )}{y_1^3},$$
$$g_{22}=\frac{2 \,x_3^2\, (3x_1^4y_2^2+6x_1^2y_1y_2+2y_1^2)}{y_1^2},\qquad \qquad g_{33}=1.$$
 The non-vanishing components $g^{ij}$ of the inverse metric tensor  are the following:
\begin{align*}
& g^{11} =\frac{\left(3 x_1^4 y_2^2+6 x_1^2 y_1 y_2+2 y_1^2\right) y_1^4}{x_3^2 x_1^2 y_2^3\left(x_1^6 y_2^3+6 x_1^4 y_1 y_2^2+12 x_1^2 y_1^2 y_2+8 y_1^3\right)} ,\\
& g^{12} =\frac{\left(2 x_1^2 y_2+3 y_1\right) y_1^3}{x_3^2 y_2\left(x_1^6 y_2^3+6 x_1^4 y_1 y_2^2+12 x_1^2 y_1^2 y_2+8 y_1^3\right)} ,\\
& g^{22} =\frac{1}{2} \frac{\left(3 x_1^2 y_2+4 y_1\right) y_1^2}{x_3^2\left(x_1^6 y_2^3+6 x_1^4 y_1 y_2^2+12 x_1^2 y_1^2 y_2+8 y_1^3\right)}, \quad g^{33} =1.
\end{align*}
Consequently, the non-vanishing components of the Cartan torsion $C_{ijk}$ are
$$\,\,C_{111} = - \frac{6\, x_1^2\, x_3^2\,   y_2^3(x_1^2 y_2+y_1)  }{y_1^5}, ~~~~~~~~~~~~~~\,\,
  C_{112}=\frac{6\, x_1^2\, x_3^2\,   y_2^2 (x_1^2 y_2+y_1)  }{y_1^4},$$
$$C_{122} =- \frac{6\, x_1^2\, x_3^2\,   y_2 (x_1^2 y_2+y_1)  }{y_1^3} , ~~~~~~~~~~~~~~
  C_{222}=\frac{6\, x_1^2\, x_3^2\,    (x_1^2 y_2+y_1)  }{y_1^2}.$$
Moreover, the geodesic spray coefficients are given by
$$G^1=\frac{(x_1 y_3-x_3 y_1)y_1}{x_1x_3}, \quad G^2=\frac{y_2y_3}{x_3},\quad  G^3=-\frac{x_3 y_2^2(x_1^4y_2^2+4x_1^2y_1y_2+4y_1^2)}{2y_1^2}.$$
 Some of the coefficients of Cartan connection are
$$\Gamma^1_{13}=\frac{1}{x_3}, \quad \Gamma^2_{23}=\frac{1}{x_3}, \quad \Gamma^3_{33}=0. $$

One can choose a concurrent $\pi$-vector field $\overline{\varphi}={\varphi }^i (x)\overline{\partial_i}=x_3 \overline{\partial_3}$ with respect to $F$, where $\overline{\partial_i}$ are the basis of fibres of $\pi^{-1}(TM)$.
 Clearly,  ${\varphi }^iC_{ijk}=0$ and
$$ {\varphi }^1_{\,\,\,|1}=\delta_1\, {\varphi }^1+{\varphi}^1\, \Gamma^1_{11}+{\varphi }^2\, \Gamma^1_{12}+{\varphi }^3\, \Gamma^1_{13}=1.$$ Similarly, ${\varphi }^2_{\,\,\,|2}=1$, ${\varphi }^3_{\,\,\,|3}=1$ and all other components of  ${\varphi }^i_{\,\,\,|j}$ vanish identically.
\\
 Moreover, the components of the corresponding $\pi$-form ${\phi}$  are
$\phi_1=\phi_2=0, \ \phi_3=x_3.$ Consequently,  $\Phi(x,y)= x_3 y_3$.
\par We then have
$$\widehat{F}(x,y)=\frac{F^2(x,y)}{F(x,y)-\Phi(x,y)}=\frac{x_3^2\left(\frac{x_1^2  y_2^2+2y_1y_2}{y_1}\right)^2+y_3^2}
{\sqrt{x_3^2\left(\frac{x_1^2  y_2^2+2y_1y_2}{y_1}\right)^2+y_3^2}- {x_3 y_3}}.$$
 Since $\widehat{F}$ satisfies \eqref{Eq:g_tilde_non_degenerate}, as it can be verified using the Maple Finsler package \cite{NF_Package}, then  $\widehat{F}$ defines a generalized $\phi$-Matsumoto metric over $D$, where  $$D:= D_{F} \cap \Big\{ \sqrt{x_3^2\left(\frac{x_1^2  y_2^2+2y_1y_2}{y_1}\right)^2+y_3^2}>{x_3 y_3}\Big\}.$$
 }
\end{example}
Let us recall the following definition.
\begin{defn}\emph{ \cite{arFinslerTaha}}
A conic-Finsler metric $F$ on  $\mathcal{D}$ is called an \textbf{almost rational} Finsler metric if all its Finsler metric tensor  components $g_{ij}(x,y)$ can be  expressed in the form
$$
g_{ij}(x,y)= \theta (x,y)\, a_{ij}(x,y),
$$
where $\theta : \mathcal{D} \longrightarrow (0,\infty )$ is a smooth function and the matrix $( a_{ij}(x,y))_{1 \leq i,j \leq n}$ is  symmetric non-degenerate and each of the functions $ a_{ij}(x,y)$ is rational  in the directional argument $y$.
\par If in addition,  $\theta $ is a rational function in $y$,   the Finsler metric $F$ is said to be a \textbf{rational} Finsler metric.
\end{defn}
\begin{rem}
The Finsler metric $F$ defined by \eqref{Lexample} is a rational Finsler metric since $g_{ij}$ can be written as $g_{ij} (x,y)= \theta  (x,y)\, a_{ij}(x,y)$ with $\theta  (x,y)= (\frac{x_3}{y_1})^2$ being a rational function in the variable $y$  and $$\qquad \qquad  \qquad a_{11}=\frac{x_1^2\, y_2^3(3x_1^2 y_2+4y_1)}{y_1^2},\quad\quad \quad a_{12}=-\frac{2 x_1^2 y_2^2 (2x_1^2\, y_2 +3\,y_1 )}{y_1^2},$$
$$a_{22}=2 (3x_1^4y_2^2+6x_1^2y_1y_2+2y_1^2),\qquad \qquad a_{33}=(\frac{y_1}{x_3})^2 $$ all are rational functions in the variable $y$.
\end{rem}
\textbf{This motivate us to  study the following problem:}
 If $g_{ij}$ are written as $g_{ij}=\theta  \, a_{ij}$ with    $a_{ij}$ are rational functions in $y$, when can $\widehat{g}_{ij}$ be written as $\widehat{g}_{ij}=\widehat{\theta }\, \widehat{a}_{ij}$ with  $ \widehat{a}_{ij}$  rational functions in $y$? Or, equivalently, under what conditions the $\phi$-Matsumoto change preserves the almost rationality property of  Finsler metrics?
\begin{lem}\label{F is a rational fn}
Let $(M,F)$ be a Finsler manifold admitting a concurrent $\pi$-vector field $\overline{\varphi}$. If $F$ is a rational function in $y$,  then $F$ is a rational Finsler metric.
\end{lem}
\begin{proof}
Suppose that the Finsler structure  $F$ is a rational function in $y$. It is clear that $F^2$ will be a rational function in $y$.  Since,  the partial differentiation $ \frac{\pa}{\pa y^{i}}$ of a geometric object that is rational in $y$ remains rational in $y$,  the Finsler metric tensor components $g_{ij}(x,y)$ are rational functions in $y$ and can be expressed in the form $g_{ij}(x,y)=\theta (x,y)\, a_{ij}(x,y)$ with  $\theta(x,y)=1 $ and all $a_{ij}(x,y)$ are rational functions in $y$.  That is,  $F$ is a rational Finsler metric.
\end{proof}
\begin{rem}
The above Lemma shows that  if $F$ is a rational function in $y$, then $F$ is a rational Finsler metric. However, the converse is not true. Indeed, \emph{Example 1} has $F$  rational Finsler metric while $F$ is not itself a rational function in $y$.
\end{rem}
\begin{thm}\label{F rational fn}
Let $(M,F)$ be a Finsler manifold admitting a concurrent $\pi$-vector field $\overline{\varphi}$.
If $F$ is a rational function in $y$,  then the  generalized $\phi$-Matsumoto metric $\widehat{F}$ is a rational Finsler metric.
\end{thm}
\begin{proof}
As $F$ is a rational function in directional argument, $F$ is a rational Finsler metric, by Lemma \ref{F is a rational fn}, that is, $g_{ij}(x,y)=\theta (x,y)\, a_{ij}(x,y)$ with  $\theta(x,y)=1$ and $a_{ij}(x,y)$ are rational functions in $y$.
  \par Under the  $\phi$-Matsumoto change \eqref{change},  the local expression of its Finsler metric components \eqref{gg3}  are given by
\begin{eqnarray}\label{gg3 locally}
    \widehat{g}_{ij}&=&  \frac{F^2(F-2\Phi)}{(F-\Phi)^3}g_{ij}+\frac{3F^4}{(F-\Phi)^4}\phi_i\,\phi_{j}
    +\frac{F^2\,\Phi(4\Phi-F)}{(F-\Phi)^4}\,\ell_{i}\,\ell_{j} \nonumber\\
    &&+\frac{F^3(F-4\Phi)}{(F-\Phi)^4}\set{\phi_{i}\,\ell_{j}+\phi_{j}\,\ell_{i}},
     \end{eqnarray}
   which can be written in the form
   \begin{eqnarray*}
  \widehat{g}_{ij}  &=& \frac{F^2}{(F-\Phi)^4}\set{(F-\Phi) (F-2\Phi) \,a_{ij} +3F^2 \phi_i\,\phi_{j}+\Phi(4\Phi-F)\,\ell_{i}\,\ell_{j}}
   \\ &&
    +\frac{F^2}{(F-\Phi)^4}  F(F-4\Phi)
    \set{\phi_{i}\,\ell_{j}+\phi_{j}\,\ell_{i}} .
    \end{eqnarray*}
\par   The following functions are obviously rational in $y$
   $$\ell_{i}\,\ell_{j} = \frac{g_{ri} y^r}{F} \frac{g_{kj} y^k}{F} =  \frac{g_{ri} y^r g_{kj} y^k}{F^2} = \frac{  a_{ri} y^r   a_{kj} y^k}{ a_{ms}y^m y^s} = \frac{ a_{ri} y^r  a_{kj} y^k}{a_{ms}y^m y^s}, $$
    $$F \ell_{i}= g_{ri} y^r =  a_{ri} y^r .$$
 \par Now,  setting
 \begin{equation}\label{tildeeta}
 \widehat{\theta } =  \frac{F^2}{(F-\Phi)^4},
  \end{equation}
 \begin{eqnarray}\label{tildeaij}
 \widehat{a}_{ij}  &= &(F^2 +2\Phi^2)  \,a_{ij} +3F^2 \phi_i\,\phi_{j}
    +4\Phi^2 \ell_{i}\,\ell_{j} -4\Phi\set{\phi_{i}\,F\,\ell_{j}+\phi_{j}\,F\,\ell_{i}}  \nonumber \\ &&
    -F\Phi(3   \,a_{ij} +\ell_{i}\,\ell_{j})+F\,\set{\phi_{i}\,F\,\ell_{j}+\phi_{j}\,F\,\ell_{i}},
 \end{eqnarray}
 imply that  $\widehat{g}_{ij} = \widehat{\theta }\,  \widehat{a}_{ij}$.  By Remark \ref{independent of the directional argument},  the functions
$\phi_{i}, \text{ for all } i=1,...,n,$ are independent of  $y$ and $\Phi$ is linear in $y$. Thereby, as $F$ is a rational function in $y$, each of the functions $ \widehat{\theta }$ and  $ \widehat{a}_{ij}$ are rational in $y$.  Hence,  $\widehat{F}$ is a rational Finsler metric.
\end{proof}
The following result provides an answer to the above mentioned question.
\begin{thm}
Let $(M,F)$ be a Finsler manifold admitting a concurrent $\pi$-vector field $\overline{\varphi}$. If $F$ is a rational Finsler metric, then the following assertions are equivalent:
\begin{description}
    \item[(a)] $F$ is a rational function in $y. \qquad\qquad$
\textbf{\emph{(b)}} $\widehat{F}$ is a rational Finsler metric.
\end{description}
\end{thm}
\begin{proof}
\textbf{(a)}$\implies$\textbf{(b)}:
It follows from Theorem \ref{F rational fn}.
\par \textbf{(b)}$\implies$\textbf{(a)}:
Now, suppose $\widehat{F}$ is a rational Finsler metric. In fact, the Finsler metric tensor $\widehat{g}_{ij}$, given by \eqref{gg3 locally}, can be written as $\widehat{g}_{ij}= \widehat{\theta }\,  \widehat{a}_{ij}$ with $\widehat{\theta } =  \frac{\theta \,F^2}{(F-\Phi)^4},$
  \begin{eqnarray*}\label{tildeaij ratL}
 \widehat{a}_{ij}  &= &  \set{ (F^2 +2\Phi^2) \,a_{ij} +3a_{rk} y^k y^r \phi_i\,\phi_{j}
    +4\Phi^2  \frac{ a_{ri} y^r  a_{kj} y^k}{a_{ms}y^m y^s}-4\Phi\set{\phi_{i}\,a_{rj} y^r +\phi_{j}\,a_{ri} y^r  }}
\\ \nonumber&&
     -F \set{ \Phi(3 a_{ij} +   \frac{ a_{ri} y^r  a_{kj} y^k}{a_{ms}y^m y^s})+\phi_{i}\,a_{rj} y^r +\phi_{j}\,a_{ri} y^r } .
 \end{eqnarray*}
 Since $F$ is a rational Finsler metric, the functions
\small{ \begin{equation}\label{aij ingeneral}
  \set{ (F^2 +2\Phi^2) \,a_{ij} +3a_{rk} y^k y^r \phi_i\,\phi_{j}
    +4\Phi^2  \frac{ a_{ri} y^r  a_{kj} y^k}{a_{ms}y^m y^s}-4\Phi\set{\phi_{i}\,a_{rj} y^r +\phi_{j}\,a_{ri} y^r  }}
 \end{equation}}
  are rational in $y$ for all $i,j$. However, the functions
 $$  F\set{ -\Phi(3 \,a_{ij} + \frac{ a_{ri} y^r  a_{kj} y^k}{a_{ms}y^m y^s})+ \phi_{i}\,a_{rj} y^r+\phi_{j}\,a_{ri} y^r} $$
are rational functions in $y$ if and only if $F$ is a rational function in $y$. That is, $ \widehat{a}_{ij} $ are rational functions in $y$ if and only if $F$ is a rational function in $y$. Which is equivalent to $ \widehat{F}$ is a rational Finsler metric if and only if $F$ is a rational function in $y$.
 \end{proof}
\begin{rem}
If $F$ is an almost rational Finsler metric, i.e.,  $\theta$ is not a rational function in $y$, we can not conclude anything about the almost rationality of $ \widehat{F}$ that can be seen from \eqref{aij ingeneral}.
\end{rem}
\textbf{Acknowledgment.} We would like to express our deep thanks to the referees for their careful reading of this manuscript and their valuable comments which led to the present version.
\providecommand{\bysame}{\leavevmode\hbox
to3em{\hrulefill}\thinspace}
\providecommand{\MR}{\relax\ifhmode\unskip\space\fi MR }
\providecommand{\MRhref}[2]{%
  \href{http://www.ams.org/mathscinet-getitem?mr=#1}{#2}
} \providecommand{\href}[2]{#2}


\begin{thebibliography}{21}
\bibitem{Erasmo} E. Caponio, M. A. Javaloyes and M.  Sánchez, 
\emph{Wind Finslerian structures: from Zermelo’s navigation to the causality of spacetimes}, Memoirs of the American Mathematical Society 300 (1501), 134.
\bibitem{r20}A.~Fr\"{o}licher and A.~Nijenhuis, \emph{Theory of vector-valued differential forms}, {I}, Ann. Proc. Kon. Ned. Akad., A,\textbf{{59}} (1956), 338--359.

\bibitem{r21}J.~Grifone, \emph{Structure pr\'esque-tangente et connexions, \textsc{I}}, Ann. Inst. Fourier, Grenoble,\textbf{ 22}, 1 (1972), 287-334.

\bibitem{h-Matsumoto change} M. K. Gupta,  A. Sahu and
 S. Sharma,  \emph{Cartan connection for h-Matsumoto change}, 2022. \url{http://arxiv.org/abs/2204.07298v1}
\bibitem{Finsler definitions} M. Javaloyes and M. Sánchez, \textit{On the definition and examples of Finsler metrics},   Ann. Sc. Norm.
Super. Pisa Cl. Sci. (5) Vol. XIII (2014), 813-858.  
\bibitem{Lovas} R.L. Lovas, \emph{A note on Finsler-Minkowski norms}, Houston J.  Math.  \textbf{33} (2007), 701-707.


\bibitem{r2.8}M.~Matsumoto, \emph{On Finsler spaces with Randers metric and special forms of important tensors}, J. Math. Kyoto Univ., \textbf{14} (1974), 477--498.

\bibitem{[6]}
 M. Matsumoto, \emph{A slope of a mountain is a Finsler surface with respect to a time measure}, J. Math. Kyoto Univ., \textbf{29} (1989) 17–25.

\bibitem{MatEguchi} M. Matsumoto and K. Eguchi, \emph{Finsler spaces admitting a concurrent vector field}, Tensor, N. S.,\textbf{ 28} (1974), 239-249.

\bibitem{Anastaise} R. Miron and M. Anastaise,  \emph{Vector bundles and Lagrange spaces with applications to Relativity}, 1997.

\bibitem{[3]}
M. Rafie-Rad and B. Rezaei, \emph{On Einstein Matsumoto metrics}, Nonlinear Anal. RWA, \textbf{13} (2012) 882–886.

\bibitem{Matsumoto Change of Finsler Metric} H.  S.  Shukla, O. P.  Pandey and H. D. Joshi, \emph{ Matsumoto change of Finsler metric},  J. Inter.  Acad. of Physical Sciences, \textbf{16}, 4 (2012), 329-341.

 \bibitem{amr2}
A. Soleiman, \emph{Energy $\beta$-conformal change in Finsler geometry}, Int. J. Geom. Meth. Mod. Phys. \textbf{9} (2012), 1250029 (21 pages).

\bibitem{Hsca} A. Soleiman and S. G. Elgendi, \emph{The geometry of Finsler spaces of Hp-scalar curvature}, Differential Geometry - Dynamical Systems, \textbf{22} (2020), 254-268.

\bibitem{Numata1} A. Soleiman and S. G. Elgendi, \emph{An intrinsic proof of Numata's theorem on Landsberg spaces},  J. Korean Math. Soc. \textbf{61 }(2024), 149-160.

\bibitem{square metric} A. Soleiman and S. G. Elgendi, \emph{On generalized Shen's square metric}, Korean J. Math. \textbf{32 }(2024), 467-484.
 \bibitem{TahaSoleiman22}A. Soleiman and E. H. Taha, \textit{Tripathi connection in Finsler geometry}, Hacet. J. Math. Stat. \textbf{51 }(2022), 136-150.

\bibitem{2014} J. Szilasi, R. L. Lovas and D. Cs Kertesz, \emph{Connections, Sprays and Finsler structures},
World Scientific, 2014.


\bibitem{Tachibana}
S. Tachibana, \emph{On Finsler spaces which admit concurrent vector field},  Tensor, N. S.  \textbf{1} (1950), 1--5.

\bibitem{Tamim-Youssef} A. A Tamim and N. L Youssef, \emph{On generalized Randers manifolds}, Algebras,
  Groups and Geometries, \textbf{16} (1999), 115-126.
\bibitem{arFinslerTaha} E. H. Taha and B. Tiwari, \emph{On almost rational Finsler metrics},  Bull.  Iran.  Math.  Soc. \textbf{49} (8) (2023).

\bibitem{Tayebi1} A. Tayebi, T. Tabatabaeifar and E. Peyghan, \emph{On the second approximate Matsumoto metric},  Bull. Korean. Math. Soc. 51(1) (2014), 115-128.
\bibitem{Tayebi2} A. Tayebi, E. Peyghan and H. Sadeghi, \emph{On Matsumoto-type Finsler metrics},
Nonlinear Analysis: RWA, 13 (2012), 2556-2561.
\bibitem{Tayebi3} A. Tayebi and M. Shahbazi Nia, \emph{On Matsumoto change of m-th root Finsler metrics},  Publications De L’institut Mathematique, tome 101(115) (2017), 183-
190.


\bibitem{r92}
N.~L. Youssef, S.~H. Abed and A.~Soleiman, \emph{Cartan and
  Berwald connections in the pullback formalism}, Algebras,
  Groups and Geometries, \textbf{25} (2008), 363-384. 

\bibitem{r99}N.~L. Youssef, S.~H. Abed and A.~Soleiman,, \emph{Intrinsic theory of projective changes in Finsler geometry}, Rend. Circ. Mat. Palermo, \textbf{60} (2011), 263–281.

\bibitem{r94a}
N.~L. Youssef, S. H. Abed and A.~Soleiman, \emph{Concurrent $\pi$-vector fields and energy
$\beta$-change}, Inter. J.  Geom. Meth.  Mod. Phys., \textbf{6} (2009) 1003–1031.

\bibitem{r94}N.~L. Youssef, S.~H. Abed and A.~Soleiman, \emph{A global approach to  the theory  of connections in Finsler geometry}, Tensor, N. S., \textbf{71} (2009),187-208. arXiv: 0801.3220 [math.DG].

\bibitem{r96}N.~L. Youssef, S.~H. Abed and A.~Soleiman, \emph{Geometric objects associated with the fundamental connections in Finsler geometry},   J. Egypt. Math. Soc.,\textbf{18} (2010), 67-90.  arXiv: 0805.2489 [math.DG].


\bibitem{NF_Package}
N.~L. Youssef and S. G. Elgendi, \emph{New Finsler   package},  Comput. Phys. Commun., \textbf{185}, 3 (2014), 986--997.


 \bibitem{Semi-concurrent}
 N. L. Youssef,   S. G. Elgendi  and E. H. Taha, \emph{Semi-concurrent vector fields in Finsler geometry},  Differ. Geom.  Appl., \textbf{65} (2019).  1--15.

\end{thebibliography}
\end{document}